\documentclass[oneside]{amsart}
\usepackage{amssymb}

\makeatletter
 \theoremstyle{plain}    
 \newtheorem{thm}{Theorem}[section]
 \numberwithin{equation}{section} 
 \numberwithin{figure}{section} 
 \theoremstyle{plain}
 \theoremstyle{plain}    
 \newtheorem*{thm*}{Theorem} 
 \theoremstyle{remark}    
 \newtheorem*{acknowledgement*}{Acknowledgement} 
 \theoremstyle{definition}
 \newtheorem{defn}[thm]{Definition}
 \theoremstyle{plain}    
 \newtheorem{cor}[thm]{Corollary} 
 \theoremstyle{plain}    
 \newtheorem{prop}[thm]{Proposition} 
 \theoremstyle{plain}    
 \newtheorem{lem}[thm]{Lemma} 

\usepackage{amsfonts}%
\makeatletter

\newcommand{\diam}{\operatorname{diam}}

\newcommand{\Ker}{\operatorname{Ker}}
\newcommand{\Lip}{\operatorname{Lip}}

\newcommand{\Aut}{\operatorname{Aut}}

\makeatother
\begin{document}

\title{Operator Algebras and Mauldin-Williams Graphs}

\author{Marius Ionescu}


\begin{abstract}
We describe a method for associating a $C^{*}$-correspondence to
a Mauldin-Williams graph and show that the Cuntz-Pimsner algebra of
this $C^{*}$-correspondence is isomorphic to the $C^{*}$-algebra
of the underlying graph. In addition, we analyze certain ideals of
these $C^{*}$-algebras.

We also investigate Mauldin-Williams graphs and fractal $C^{*}$-algebras
in the context of the Rieffel metric. This generalizes the work of
Pinzari, Watatani and Yonetani. Our main result here is a {}``no
go'' theorem showing that such algebras must come from the commutative
setting.
\end{abstract}
\maketitle

\section{Introduction}

In recent years many classes of $C^{\ast}$-algebras have been shown
to fit into the Pimsner construction of what are known now as Cuntz-Pimsner
algebras (see \cite{Pimnser} and \cite{Muhly-Baruc-tensoalgebras}).
His construction is based on a so-called $C^{\ast}$-correspondence
over a $C^{\ast}$-algebra. For example, a natural $C^{\ast}$-correspondence
can be associated with a graph $G$ (see \cite{Fowler-Laca-Raeburn}
and \cite[Example 1.5]{Fowler-Muhly-Raeburn}). The Cuntz-Pimsner
algebra of this $C^{\ast}$-correspondence is isomorphic to the graph
$C^{\ast}$-algebra $C^{\ast}(G)$ as defined in \cite{Kumjian-Pask-Raeburn-directedgraphs}.
Another example is the $C^{\ast}$-correspondence associated with
a local homeomorphism on a compact metric space studied by Deaconu
in \cite{Deaconu2}, and the $C^{\ast}$-correspondence associated
with a local homeomorphism on a locally compact space studied by Deaconu,
Kumjian, and Muhly in \cite{Deaconu-Kumjian-Muhly}. They showed that
the Cuntz-Pimsner algebra is isomorphic to the groupoid $C^{\ast}$-algebra
associated with a local homeomorphism in \cite{Deaconu1}, \cite{Deaconu-Kumjian-Muhly}
and \cite{Renault}.

By a (directed) \emph{graph} we mean a system $G=(V,E,r,s)$ where
$V$ and $E$ are finite sets, called the sets of \emph{vertices}
and \emph{edges} (respectively) of the graph and where $r$ and $s$
are maps from $E$ to $V$, called the \emph{range} and \emph{source}
maps, respectively. Thus, $s(e)$ is the source of an edge $e$ and
$r(e)$ is its range. A \emph{Mauldin-Williams graph} is a graph $G$
together with a collection of compact metric spaces, one for each
vertex of the graph, and a collection of contraction maps, one for
each edge of the graph which satisfy certain properties (see Definition
\ref{def:Mauldin-Williams-Graphs} below). In this note we follow
the notations from \cite{Edgar}. We associate with such a system
a $C^{\ast}$-correspondence which reflects the dynamics of the Mauldin-Williams
graph and we analyze the Cuntz-Pimsner algebra of this $C^{\ast}$-correspondence.
Our construction is related with topological generalizations of graph
$C^{\ast}$-algebras of Katsura \cite{Ka} and Muhly and Tomforde
\cite{MT}. The study of the Cuntz-Pimsner algebra associated with
graph dynamical systems was initiated in \cite{Pinzari-Watatami-Y},
where the authors consider the case when the graph $G$ consists of
a single vertex $v$ and edges $e_{1},e_{2},\dots,e_{n}$. In this
case, the $\phi_{e}$'s constitute what is known as in \emph{iterated
function system} acting on the space $(T_{v},\rho_{v})$. They conclude
that $\mathcal{O}(\mathcal{X})$ is isomorphic to the Cuntz algebra
$\mathcal{O}_{n}$. The first of our two principal theorems in this
note generalizes this result. Our proof is different from the proof
in \cite{Pinzari-Watatami-Y} and reveals extra structure.

\begin{thm*}
Let $(G,\{ T_{v},\rho_{v}\}_{v\in V},\{\phi_{e}\}_{e\in E})$ be a
Mauldin-Williams graph such that the graph $G$ has no sinks and no
sources. Let $A$ and $\mathcal{X}$ be the associated $C^{\ast}$-algebra
and $C^{\ast}$-correspondence. Then the Cuntz-Pimsner algebra $\mathcal{O}(\mathcal{X})$
is isomorphic to $C^{\ast}(G)$ of \cite{Cuntz-Krieger}. 
\end{thm*}
Thus the structure of $\mathcal{O}(\mathcal{X})$ is completely determined
by the graph $G$. From one perspective, this result is somewhat disappointing.
Given the richness of dynamical systems expressed as Mauldin-Williams
graphs and given the fact that Cuntz-Pimsner algebras generalize crossed
products, one might expect a lively interplay between the dynamics
and the structure of $\mathcal{O}(\mathcal{X})$. However, the {}``rigidity''
that this theorem expresses is quite remarkable and it may inspire
one to wonder about the natural limits of the result.

In particular, one might wonder if there are noncommutative versions
of Mauldin-Williams graphs and whether these might prove to have a
richer theory. This thought was taken up in \cite[Section 4.3]{Pinzari-Watatami-Y}
where Pinzari, Watatani and Yonetani considered noncommutative iterated
function systems based on Rieffel's notion of {}``noncommutative
metric spaces'' \cite{Rieffel1,Rieffel2}. The second objective of
this note is to show that noncommutative iterated function systems
of Pinzari, Watatani and Yonetani can be formulated in the setting
of Mauldin-Williams-type graphs but that the generality gained is
illusory. Roughly, the Rieffel metric is a metric on the state space
of a (not necessarily commutative) $C^{\ast}$-algebra $A$ defined
by a certain subset of {}``Lipschitz elements'' in $A$ (see Definition
\ref{def:(The-Rieffel-metric)} below). When we associate a $C^{\ast}$-algebra
$A_{v}$ to each vertex $v\in V$, for a prescribed graph $G=(V,E,r,s)$,
and when we associate a $\ast$-homomorphism $\phi_{e}:A_{s(e)}\rightarrow A_{r(e)}$
to each edge $e\in E$ that is strictly contractive with respect to
the Rieffel metrics on $A_{s(e)}$ and $A_{r(e)}$, we call the resulting
system a noncommutative Mauldin-Williams graph. It also gives rise
to a Cuntz-Pimsner algebra $\mathcal{O}(\mathcal{X})$. Our second
objective in this note is to show that once more $\mathcal{O}(\mathcal{X})$
is isomorphic to the Cuntz-Krieger algebra associated to $G$. In
fact, we shall show in Theorem \ref{thm:homintocomm} that in such
situations, the $C^{\ast}$-algebras $A_{v}$ are necessarily commutative.
This implies, in particular, that the structures considered by Pinzari,
Watatani and Yonetani are necessarily no more general than those arising
from ordinary iterated function systems.

\begin{acknowledgement*}
The author would like to thank Mark Tomforde for many useful comments
regarding the work on this paper and especially Paul Muhly for his
constant encouragement and beneficial conversations. The author would
also like to express his thanks to the referee for pointing out errors
in the original draft of this note and for suggesting the use of the
Gauge-Invariant Uniqueness Theorem in the proof of Theorem \ref{thm:CuntzKrieger},
and to Yasuo Watatani, Alex Kumjian, and Dorin Dutkay for useful conversations
after the first version of the paper was written.
\end{acknowledgement*}

\section{The Cuntz-Pimsner algebra associated to a Mauldin-Williams graph}

\begin{defn}
\label{def:Mauldin-Williams-Graphs}By a \emph{Mauldin-Williams} \emph{graph}
(see \cite{Mauldin-Williams} and \cite{Edgar}), we mean a system
$(G,\{ T_{v},\rho_{v}\}_{v\in V},\{\phi_{e}\}_{e\in E})$ where $G=(V,E,r,s)$
is a graph and where $\{ T_{v},\rho_{v}\}_{v\in V}$ and $\{\phi_{e}\}_{e\in E}$
are families such that:
\begin{enumerate}
\item For each $v\in V$, $T_{v}$ is a compact metric space with a prescribed
metric $\rho_{v}$.
\item For $e\in E$, $\phi_{e}$ is a continuous map from $T_{r(e)}$ to
$T_{s(e)}$ such that\[
\rho_{s(e)}(\phi_{e}(x),\phi_{e}(y))\leq c\rho_{r(e)}(x,y)\]
 for some constant $c$ satisfying $0<c<1$ (independent of $e$)
and all $x,y\in T_{r(e)}$.
\end{enumerate}
We shall assume, too, that the functions $s$ and $r$ are surjective.
Thus, we assume that there are no sinks and no sources in the graph
$G$. An \emph{invariant list} associated with a Mauldin-Williams
graph is a family $(K_{v})_{v\in V}$ of compact sets, such that $K_{v}\subset T_{v}$
for all $v\in V$ and\begin{align*}
K_{v} & =\bigcup_{e\in E,s(e)=v}\phi_{e}(K_{r(e)}).\end{align*}
Since each $\phi_{e}$ is a contraction, a Mauldin-Williams graph
$(G,\{ T_{v},\rho_{v}\}_{v\in V},\{\phi_{e}\}_{e\in E})$ has a unique
invariant list (see \cite[Theorem 1]{Mauldin-Williams}). We set $T:=\bigcup_{v\in V}T_{v}$
and $K:=\bigcup_{v\in V}K_{v}$ and we call $K$ the \emph{invariant
set} of the Mauldin-Williams graph.
\end{defn}
In the particular case when we have one vertex $v$ and $n$ edges,
i.e. in the setting of an \emph{iterated function system,} the invariant
set is the unique compact subset $K:=K_{v}$ of $T=T_{v}$ such that\[
K=\phi_{1}(K)\cup\cdots\cup\phi_{n}(K).\]

\begin{defn}
\label{def:Mauldinwilliams module}Given a Mauldin-Williams graph
$(G,\{ T_{v},\rho_{v}\}_{v\in V},\{\phi_{e}\}_{e\in E})$ we construct
a so-called $C^{*}$-correspondence $\mathcal{X}$ over the $C^{*}$-algebra
$A=C(T)$, where $T=\coprod_{v\in V}T_{v}$ is the disjoint union
of the spaces $T_{v}$, as follows. Let $E\times_{G}T=\{(e,x)\:|\: x\in T_{r(e)}\}$.
Then, by our finiteness assumptions, $E\times_{G}T$ is a compact
space. We set $\mathcal{X}=C(E\times_{G}T)$ and view $\mathcal{X}$
as a bimodule over $C(T)$ via the formulae:\[
\xi\cdot a(e,x):=\xi(e,x)a(x)\]
 and\[
a\cdot\xi(e,x):=a\circ\phi_{e}(x)\xi(e,x),\]
 where $a\in C(T)$ and $\xi\in C(E\times_{G}T)$. Further, $\mathcal{X}$
comes equipped with the structure of a Hilbert $C^{*}$-module over
$C(T)$ via the formula\[
\langle\xi,\eta\rangle_{A}(x):=\sum_{{\genfrac{}{}{0pt}{}{e\in E}{x\in T_{r(e)}}}}\overline{\xi(e,x)}\eta(e,x)\]
 for all $\xi,\eta\in\mathcal{X}$, so that, in the language of \cite{Muhly-Baruc-tensoalgebras},
$\mathcal{X}$ may be viewed as a $C^{*}$-\emph{correspondence over}
$C(T)$. Since there are no sources in the graph $G$, the $A$-valued
inner product is well defined. Let $n=|E|$ and let $C^{n}(A)$ be
the column space over $A$, i.e. $C^{n}(A)=\{(\xi_{e})_{e\in E}\;:\;\xi_{e}\in A,\mbox{for all}\: e\in E\}$.
Then we view $\mathcal{X}$ as a subset of $C^{n}(A)$.

We note that the left action is given by the $\ast$-homomorphism
$\Phi:A\rightarrow\mathcal{L}(\mathcal{X})$, $(\Phi(a)\xi)(e,x)=a\circ\phi_{e}(x)\xi(e,x)$.
Then $\Phi$ is faithful if and only if $K=T$. 
\end{defn}
In the case of an iterated function system the $C^{\ast}$-correspondence
will be the full column space over the $C^{\ast}$-algebra $A=C(T)$,
that is $\mathcal{X}=C^{n}(A)$, with the structure:

\begin{itemize}
\item The right action is the untwisted right multiplication, i.e. $\xi\cdot a(i,x)=\xi(i,x)a(x)$
for all $i\in\{1,\dots,n\}$ and $x\in T$.
\item The left action is given by the $\ast$-homomorphism $\Phi:A\rightarrow\mathcal{L}(\mathcal{X})$
defined by the formula $\Phi(a)(\xi)(i,x)=a\circ\varphi_{i}(x)\xi_{i}(x)$.
\item The $A$-valued inner product given by the formula:\[
\langle\xi,\eta\rangle_{A}(x)=\sum_{i=1}^{n}\xi^{*}(i,x)\eta(i,x).\]

\end{itemize}
Given a $C^{*}$-correspondence $\mathcal{X}$ over a $C^{*}$-algebra
$A$ a \emph{Toeplitz representation} of $\mathcal{X}$ in a $C^{*}$-algebra
$B$ consists of a pair $(\psi,\pi)$, where $\psi:\mathcal{X}\rightarrow B$
is a linear map and $\pi:A\rightarrow B$ is a $\ast$-homomorphism
such that\[
\psi(x\cdot a)=\psi(x)\pi(a)\;,\;\psi(a\cdot x)=\pi(a)\psi(x),\]
 i.e. the pair $(\psi,\pi)$ is a bimodule map and\[
\psi(x)^{*}\psi(y)=\pi(\langle x,y\rangle_{A}).\]
 That is, the map $\psi$ preserves inner product (see \emph{\cite[Section 1]{Fowler-Muhly-Raeburn}}).
Given such a Toeplitz representation, there is a $\ast$-homomorphism
$\pi^{(1)}$ from $\mathcal{K}(\mathcal{X})$ into $B$ which satisfies\[
\pi^{(1)}(\Theta_{x,y})=\psi(x)\psi(y)^{*}\;\;\mbox{for}\;\mbox{all}\; x,y\in\mathcal{X},\]
 where $\Theta_{x,y}=x\otimes\tilde{y}$ is the rank one operator
defined by $\Theta_{x,y}(z)=x\cdot\langle y,z\rangle_{A}$.

We define then \[
J(\mathcal{X}):=\Phi^{-1}(\mathcal{K}(\mathcal{X})),\]
 which is a closed two sided-ideal in $A$ (see \cite[Definition 1.1]{Fowler-Muhly-Raeburn}).
Let $K$ be an ideal in $J(\mathcal{X})$. We say that a Toeplitz
representation $(\psi,\pi)$ of $\mathcal{X}$ is \emph{coisometric}
on $K$ if\[
\pi^{(1)}(\Phi(a))=\pi(a)\;\mbox{for}\;\mbox{all}\; a\in K.\]
 When $(\psi,\pi)$ is coisometric on all of $J(\mathcal{X})$, we
say that it is \emph{Cuntz-Pimsner covariant.}

It is shown in \cite[Proposition 1.3]{Fowler-Muhly-Raeburn} that,
for an ideal $K$ in $J(\mathcal{X})$, there is a $C^{*}$-algebra
$\mathcal{O}(K,\mathcal{X})$ and a Toeplitz representation $(k_{\mathcal{X}},k_{A})$
of $\mathcal{X}$ into $\mathcal{O}(K,\mathcal{X})$ which is coisometric
on \emph{$K$} and satisfies:

\begin{enumerate}
\item for every Toeplitz representation $(\psi,\pi)$ of $\mathcal{X}$
which is coisometric on $K$, there is a $\ast$-homomorphism $\psi\times_{K}\pi$
of $\mathcal{O}(K,\mathcal{X})$ such that $(\psi\times_{K}\pi)\circ k_{\mathcal{X}}=\psi$
and $(\psi\times_{K}\pi)\circ k_{A}=\pi$; and
\item $\mathcal{O}(K,\mathcal{X})$ is generated as a $C^{*}$-algebra by
$k_{\mathcal{X}}(\mathcal{X})\cup k_{A}(A)$. 
\end{enumerate}
The algebra $\mathcal{O}(\{0\},\mathcal{X})$ is the \emph{Toeplitz
algebra} $\mathcal{T}_{\mathcal{X}}$, and $\mathcal{O}(J(\mathcal{X}),\mathcal{X})$
is the \emph{Cuntz-Pimsner algebra $\mathcal{O}_{\mathcal{X}}$.}

For a finite graph $G=(E,V,r,s)$, a \emph{Cuntz-Krieger $G$-family}
consists of a family $\{ P_{v}\;:v\in V\}$ of mutually orthogonal
projections and a family of partial isometries $\{ S_{e}\}_{e\in E}$
such that\[
S_{e}^{\ast}S_{e}=P_{r(e)}\;\mbox{for }e\in E,\mbox{ and }P_{v}=\sum_{s(f)=v}S_{f}S_{f}^{\ast}\mbox{ for }v\in s(E).\]
 The edge matrix of $G$ is the $E\times E$ matrix defined by\[
A_{G}(e,f)=\begin{cases}
1\;\mbox{if}\; r(e)=s(f)\\
0\;\mbox{otherwise}.\end{cases}\]
 Then, a Cuntz-Krieger $G$-family satisfies:\[
S_{e}^{\ast}S_{e}=\sum_{f\in E}A_{G}(e,f)S_{f}S_{f}^{\ast}\]
 for every $e\in E$ such that $A_{G}(e,\cdot)$ has nonzero entries.
It is shown in \cite[Theorem 1.2]{Kumjian-Pask-Raeburn-directedgraphs}
that there exists a $C^{\ast}$-algebra $C^{\ast}(G)$ generated by
a Cuntz-Krieger $G$-family $\{ S_{e},P_{v}\}$ of non-zero elements
such that, for every Cuntz-Krieger $G$-family $\{ W_{e},T_{v}\}$
of partial isometries on $H$, there is a representation $\pi$ of
$C^{\ast}(G)$ on $H$ such that $\pi(S_{e})=W_{e}$ and $\pi(P_{v})=T_{v}$
for all $e\in E$ and $v\in V$. The triple $(C^{\ast}(G),S_{e},P_{v})$
is unique up to isomorphism. Since we are assuming that $G$ has no
sinks, $\{ S_{e}\}_{e\in E}$ is a Cuntz-Krieger family for the edge
matrix $A_{G}$ in the sense of \cite{Cuntz-Krieger} (see \cite[Section 1]{Kumjian-Pask-Raeburn-directedgraphs}),
and the projections $P_{v}$ are redundant. If the matrix $A_{G}$
satisfies Condition (I) from \cite{Cuntz-Krieger} (or, equivalently,
since $G$ is finite, if $G$ satisfies Condition (L) from \cite{Kumjian-Pask-Raeburn-directedgraphs},
which asserts that every loop has an exit), then $C^{\ast}(G)$ is
unique and is isomorphic to the Cuntz-Krieger algebra from \cite{Cuntz-Krieger}
(see \cite[Theorem 3.7]{Kumjian-Pask-Raeburn-directedgraphs}).

\begin{thm}
\label{thm:CuntzKrieger}Let $(G,\{ T_{v},\rho_{v}\}_{v\in V},\{\phi_{e}\}_{e\in E})$
be a Mauldin-Williams graph such that the graph $G$ has no sinks
and no sources. Let $A$ and $\mathcal{X}$ be defined as above. Then
the Cuntz-Pimsner algebra $\mathcal{O}_{\mathcal{X}}$ is isomorphic
to $C^{\ast}(G)$. 
\end{thm}
Before proving the theorem, we introduce some notation. For $k\geq2$,
set\[
E^{k}:=\{\alpha=(\alpha_{1},\cdots,\alpha_{k})\;:\;\alpha_{i}\in E\;\mbox{and}\; r(\alpha_{i})=s(\alpha_{i+1}),i=1,\cdots,k-1\},\]
 the set of \emph{paths of length $k$} in the graph $G$. Let $E^{\ast}=\bigcup_{k\in N}E^{k}$,
the space of \emph{finite paths} in the graph $G$. Also the \emph{infinite
path space $E^{\infty}$} is defined to be \[
E^{\infty}:=\{(\alpha_{i})_{i\in\mathbb{N}}\;:\;\alpha_{i}\in E\;\mbox{and}\; r(\alpha_{i})=s(\alpha_{i+1})\;\mbox{for all}\: i\in\mathbb{N}\}.\]
For $v\in V$, we also define $E^{k}(v):=\{\alpha\in E^{k}\;:\; s(\alpha)=v\}$,
and we define $E^{\ast}(v)$ and $E^{\infty}(v)$ in a similar way.
We consider $E^{\infty}(v)$ endowed with the metric $\delta_{v}(\alpha,\beta)=c^{|\alpha\wedge\beta|}$
if $\alpha\ne\beta$ and $0$ otherwise, where $\alpha\wedge\beta$
is the longest common prefix of $\alpha$ and $\beta$, and $|w|$
is the length of the word $w\in E^{\ast}$ (see \cite[Page 116]{Edgar}).
Then $E^{\infty}(v)$ is a compact metric space, and, since $E^{\infty}$
equals the disjoint union of the spaces $E^{\infty}(v)$, $E^{\infty}$
becomes a compact metric space in a natural way. 

For $\alpha\in E^{k}$, we write $\phi_{\alpha}=\phi_{\alpha_{1}}\circ\cdots\circ\phi_{\alpha_{k}}$
and $S_{\alpha}=S_{\alpha_{1}}\cdots S_{\alpha_{k}}$. Let $\mathcal{S}_{v}$
be the state space of $A_{v}=C(T_{v})$ and $\mathcal{S}=\prod_{v\in V}\mathcal{S}_{v}$.
We consider the metrics $L_{v}$ defined on $\mathcal{S}_{v}$ by
the formula\begin{equation}
L_{v}(\mu,\nu)=\sup\{|\mu(f)-\nu(f)|\quad:\quad f\in\operatorname{Lip}(T_{v}),c_{f}\leq1\},\label{eq:HucthmetricMW}\end{equation}
 where $\operatorname{Lip}(T_{v})$ is the space of Lipschitz functions
on $T_{v}$ and $c_{f}$ is the Lipschitz constant of the Lipschitz
function $f$. For $f\in\operatorname{Lip}(T_{v})$, $v\in V$ and
$\mu,\nu\in\mathcal{S}_{v}$\begin{equation}
|\mu(f)-\nu(f)|\leq c_{f}L_{v}(\mu,\nu).\label{eq:lema1CK}\end{equation}
 Further, if $\alpha\in E^{k}$, $k\geq1$, and $\mu,\nu\in\mathcal{S}_{r(\alpha)}$,
then $\mu\circ\phi_{\alpha_{k}}^{-1}\circ\cdots\circ\phi_{\alpha_{1}}^{-1},\nu\circ\phi_{\alpha_{k}}^{-1}\circ\cdots\circ\phi_{\alpha_{1}}^{-1}\in\mathcal{S}_{s(\alpha)}$
and\begin{equation}
L_{s(\alpha_{1})}(\mu\circ\phi_{\alpha_{k}}^{-1}\circ\cdots\circ\phi_{\alpha_{1}}^{-1},\nu\circ\phi_{\alpha_{k}}^{-1}\circ\cdots\circ\phi_{\alpha_{1}}^{-1})\leq c^{k}\operatorname{diam}_{L}(\mathcal{S}),\label{eq:lemmaCK}\end{equation}
 where $\operatorname{diam}_{L}(\mathcal{S})=\max_{v\in V}\operatorname{diam}_{L_{v}}(\mathcal{S}_{v})$.

For $\alpha\in E^{\infty}$, the sequence $\left(\phi_{\alpha_{1}\dots\alpha_{k}}(T_{r(\alpha_{k})})\right)_{k\in\mathbb{N}}\subset T_{s(\alpha)}$
is a decreasing sequence of compact sets. Moreover, $\diam\left(\phi_{\alpha_{1}\dots\alpha_{k}}(T_{r(\alpha_{k})})\right)\le c^{k}D$,
where $D:=\max_{v\in V}\diam(T_{v})$. Therefore $\lim_{k\to\infty}\diam\left(\phi_{\alpha_{1}\dots\alpha_{k}}(T_{r(\alpha_{k})})\right)=0$,
so the intersection $\bigcap_{k\in\mathbb{N}}\phi_{\alpha_{1}\dots\alpha_{k}}(T_{r(\alpha_{k})})$
consists of a single point, $x_{\alpha}\in T_{s(\alpha)}$. Hence
we can define a map $\Pi:E^{\infty}\rightarrow T$ by the formula\[
\Pi(\alpha)=x_{\alpha}.\]
Then $\Pi$ is a continuous map and its image is the invariant set
$K$ of the Mauldin-Williams graph.

\noindent \emph{Proof} \emph{of the Theorem \ref{thm:CuntzKrieger}}.

Let $\mu_{0}=(\mu_{v}^{0})_{v\in V}\in\mathcal{S}$ be fixed and $a=\sum_{v\in V}^{\oplus}a_{v}\in\operatorname{Lip}(T)$.
We define\[
i_{A}(a)=\lim_{k\rightarrow\infty}\sum_{\alpha\in E^{k}}\mu_{r(\alpha)}^{0}(a_{s(\alpha)}\circ\phi_{\alpha})S_{\alpha}S_{\alpha}^{*}.\]
We prove that $i_{A}$ is a norm decreasing $\ast$-homomorphism from
the $\ast$-algebra $\Lip(T)$ into $C^{\ast}(G)$. Then, since $\Lip(T)$
is a dense $\ast$-subalgebra of $A=C(T)$, we can extend $i_{A}$
to $A$.

We show first that the limit from the definition of $i_{A}(a)$ exists.
Let $a\in\operatorname{Lip}(T)$ and let $\varepsilon>0$. Choose
$k_{0}\in\mathbb{N}$ such that $c^{k}\operatorname{diam}_{L}(\mathcal{S})c_{a}<\varepsilon\quad\mbox{for all}\: k\geq k_{0}$.
Set $a_{k}:=\sum_{\alpha\in E^{k}}\mu_{r(\alpha)}^{0}(a_{s(\alpha)}\circ\phi_{\alpha})S_{\alpha}S_{\alpha}^{*}$.
Let $k,m\geq k_{0}$ and suppose that $k>m$. Then\begin{align*}
a_{m}-a_{k}= & \sum_{\alpha\in E^{m}}\sum_{{\genfrac{}{}{0pt}{}{\beta\in E^{k-m}}{s(\beta)=r(\alpha)}}}\mu_{r(\alpha)}^{0}(a_{s(\alpha)}\circ\phi_{\alpha})S_{\alpha}S_{\beta}S_{\beta}^{*}S_{\alpha}^{*}\\
- & \sum_{\alpha\in E^{m}}\sum_{{\genfrac{}{}{0pt}{}{\beta\in E^{k-m}}{s(\beta)=r(\alpha)}}}\mu_{r(\beta)}^{0}(a_{s(\alpha)}\circ\phi_{\alpha\beta})S_{\alpha}S_{\beta}S_{\beta}^{*}S_{\alpha}^{*}\\
= & \sum_{\alpha\in E^{m}}\sum_{{\genfrac{}{}{0pt}{}{\beta\in E^{k-m}}{s(\beta)=r(\alpha)}}}(\mu_{r(\alpha)}^{0}(a_{s(\alpha)}\circ\phi_{\alpha})-\mu_{r(\beta)}^{0}(a_{s(\alpha)}\circ\phi_{\alpha\beta}))S_{\alpha}S_{\beta}S_{\beta}^{*}S_{\alpha}^{*}\\
= & \sum_{\alpha\in E^{m}}\sum_{{\genfrac{}{}{0pt}{}{\beta\in E^{k-m}}{s(\beta)=r(\alpha)}}}(\mu_{r(\alpha)}^{0}(a_{s(\alpha)}\circ\phi_{\alpha})-\mu_{r(\beta)}^{0}\circ\phi_{\beta}^{-1}(a_{s(\alpha)}\circ\phi_{\alpha}))S_{\alpha}S_{\beta}S_{\beta}^{*}S_{\alpha}^{*}.\end{align*}
Since $|\mu_{r(\alpha)}^{0}(a_{s(\alpha)}\circ\phi_{\alpha})-\mu_{r(\beta)}^{0}\circ\phi_{\beta}^{-1}(a_{s(\alpha)}\circ\phi_{\alpha})|<\varepsilon$,
by Equations (\ref{eq:lemmaCK}) and (\ref{eq:lema1CK}), for all
$\alpha\in E^{m},\beta\in E^{k-m}$ such that $s(\beta)=r(\alpha)$,
and since $S_{\alpha}S_{\beta}S_{\beta}^{*}S_{\alpha}^{*}$ are orthogonal
projections, $\left\Vert a_{m}-a_{k}\right\Vert <\varepsilon$, for
all $m,k\geq k_{0}$. So $(a_{k})_{k\in\mathbb{N}}$ is a Cauchy sequence,
hence convergent. Since $\Vert a_{m}\Vert=\max_{\alpha\in E^{m}}|\mu_{r(\alpha)}^{0}(a_{s(\alpha)}\circ\phi_{\alpha})|\le\Vert a\Vert$
for all $m\in\mathbb{N}$, $\Vert i_{A}(a)\Vert\le\Vert a\Vert$ for
all $a\in A$.

Next we prove that $i_{A}$ is a homomorphism. Let $a,b\in\operatorname{Lip}(T)$.
Then for each $\alpha\in E^{\infty}$ there is a point $x_{\alpha}\in K$
such that $\bigcap_{k\in\mathbb{N}}\phi_{\alpha_{1}\dots\alpha_{k}}(T_{r(\alpha_{k})})=\{ x_{\alpha}\}$.
Then $\lim_{k\rightarrow\infty}\mu_{r(\alpha_{k})}^{0}(a_{s(\alpha)}\circ\phi_{\alpha_{1}\dots\alpha_{k}})=a(x_{\alpha})$,
$\lim_{k\rightarrow\infty}\mu_{r(\alpha_{k})}^{0}(b_{s(\alpha)}\circ\phi_{\alpha_{1}\dots\alpha_{k}})=b(x_{\alpha})$
and $\lim_{k\rightarrow\infty}\mu_{r(\alpha_{k})}^{0}((ab)_{s(\alpha)}\circ\phi_{\alpha_{1}\dots\alpha_{k}})=a(x_{\alpha})b(x_{\alpha})$.
Let $\varepsilon>0$. Since $\diam(\phi_{\alpha_{1}\dots\alpha_{k}}(T_{r(\alpha_{k})}))\le c^{k}D$
for all $\alpha\in E^{\infty}$ and $k\in\mathbb{N}$, there exists
some $N\in\mathbb{N}$ such that $|\mu_{r(\alpha)}(a_{s(\alpha)}\circ\phi_{\alpha_{1}\dots\alpha_{k}})-a(x_{\alpha})|<\varepsilon$,
$|\mu_{r(\alpha)}(b_{s(\alpha)}\circ\phi_{\alpha_{1}\dots\alpha_{k}})-b(x_{\alpha})|<\varepsilon$
and $|\mu_{r(\alpha)}(a_{s(\alpha)}\circ\phi_{\alpha_{1}\dots\alpha_{k}}b\circ\phi_{\alpha_{1}\dots\alpha_{k}})-a(x_{\alpha})b(x_{\alpha})|<\varepsilon$
for all $k\ge N$ and \emph{for all} $\alpha\in E^{\infty}$. Then\[
\left\Vert \sum_{\alpha\in E^{k}}\mu_{r(\alpha)}^{0}(a_{s(\alpha)}\circ\phi_{\alpha}b_{s(\alpha)}\circ\phi_{\alpha})S_{\alpha}S_{\alpha}^{\ast}-\sum_{\alpha\in E^{k}}\mu_{r(\alpha)}^{0}(a_{s(\alpha)}\circ\phi_{\alpha})\mu_{r(\alpha)}^{0}(b_{s(\alpha)}\circ\phi_{\alpha})S_{\alpha}S_{\alpha}^{\ast}\right\Vert \]
\[
\le\max_{\alpha\in E^{k}}|\mu_{r(\alpha)}^{0}(a_{s(\alpha)}\circ\phi_{\alpha}b_{s(\alpha)}\circ\phi_{\alpha})-\mu_{r(\alpha)}^{0}(a_{s(\alpha)}\circ\phi_{\alpha})\mu_{r(\alpha)}^{0}(b_{s(\alpha)}\circ\phi_{\alpha})|<\varepsilon(1+\Vert a\Vert+\Vert b\Vert)\]
for all $k\ge N$. Thus $i_{A}(ab)=i_{A}(a)i_{A}(b)$. Hence $i_{A}$
is a homomorphism and one can easily see that it is a $\ast$-homomorphism.

Further $i_{A}$ satisfies the equation\begin{equation}
i_{A}(a)S_{e}=S_{e}i_{A}(a_{s(e)}\circ\phi_{e})\quad\mbox{for all }a\in A,\mbox{ and }e\in E,\label{eq:covarrel}\end{equation}
 where we extend the map $a_{s(e)}\circ\phi_{e}$ to all $T$ by setting
it to be $0$ when $x\notin T_{r(e)}$, since, for $a\in\operatorname{Lip}(T)$,
we have:\begin{align*}
i_{A}(a)S_{e}= & \left(\lim_{k\rightarrow\infty}\sum_{\alpha\in E^{k}}\mu_{r(\alpha)}^{0}(a_{s(\alpha)}\circ\phi_{\alpha})S_{\alpha}S_{\alpha}^{*}\right)S_{e}\\
= & \lim_{k\rightarrow\infty}\sum_{\alpha\in E^{k},\alpha_{1}=e}\mu_{r(\alpha)}^{0}(a_{s(\alpha)}\circ\phi_{\alpha})S_{e}S_{\alpha_{2}}\cdots S_{\alpha_{k}}S_{\alpha_{k}}^{*}\cdots S_{\alpha_{2}}^{*}\\
= & S_{e}\lim_{k\rightarrow\infty}\sum_{\alpha^{\prime}\in E^{k-1}(r(e))}\mu_{r(\alpha^{\prime})}^{0}(a_{s(e)}\circ\phi_{e}\circ\phi_{\alpha^{\prime}})S_{\alpha^{\prime}}S_{\alpha^{\prime}}^{*}=S_{e}i_{A}(a_{s(e)}\circ\phi_{e}).\end{align*}
 We also define the linear map $i_{\mathcal{X}}:\mathcal{X}\rightarrow C^{*}(G)$
by the formula\[
i_{\mathcal{X}}(\xi)=\sum_{e\in E}S_{e}i_{A}(\xi_{e}),\]
 where $\xi_{e}\in C(T)$ is defined by $\xi_{e}(x)=\xi(e,x)$ if
$x\in T_{r(e)}$ and $0$ otherwise. We have\begin{align*}
i_{\mathcal{X}}(\xi\cdot a)= & \sum_{e\in E}S_{e}i_{A}(\xi_{e}a)=\sum_{e\in E}S_{e}i_{A}(\xi_{e})i_{A}(a)=i_{\mathcal{X}}(\xi)i_{A}(a),\end{align*}
\begin{align*}
i_{\mathcal{X}}(a\cdot\xi)= & \sum_{e\in E}S_{e}i_{A}(a_{s(e)}\circ\phi_{e}\xi_{e})=\sum_{e\in E}S_{e}i_{A}(a_{s(e)}\circ\phi_{e})i_{A}(\xi_{e})\\
= & \sum_{e\in E}i_{A}(a)S_{e}i_{A}(\xi_{e})=i_{A}(a)i_{\mathcal{X}}(\xi)\end{align*}
and\begin{align*}
i_{\mathcal{X}}(\xi)^{*}i_{\mathcal{X}}(\eta)= & \left(\sum_{e\in E}S_{e}i_{A}(\xi_{e})\right)^{*}\left(\sum_{f\in E}S_{f}i_{A}(\eta_{f})\right)\\
= & \sum_{e\in E}i_{A}(\xi_{e})^{*}i_{A}(\eta_{e})=i_{A}\left(\sum_{e\in E}\xi_{e}^{*}\eta_{e}\right)=i_{A}(\langle\xi,\eta\rangle_{A}).\end{align*}
 Hence $(i_{A},i_{\mathcal{X}})$ is a Toeplitz representation.

Let $J(\mathcal{X}):=\Phi^{-1}(\mathcal{K}(\mathcal{X}))$. Note that
$J(\mathcal{X})=A$ since for $a\in A$ we have\[
\Phi(a)\xi=\sum_{e\in E}\Theta_{x^{e},\delta^{e}}(\xi),\]
 where $x^{e}\in\mathcal{X}$ is defined by $x_{f}^{e}=a_{s(e)}\circ\phi_{e}\delta_{f}^{e}$,\[
\delta^{e}(f,x)=\begin{cases}
1\;\mbox{if}\quad f=e\\
0\;\mbox{otherwise}.\end{cases}\]
Then, for $a\in A$, we have\begin{align*}
i_{A}^{(1)}(\Phi(a))= & i_{A}^{(1)}\left(\sum_{e\in E}\Theta_{x^{e},\delta^{e}}\right)=\sum_{e\in E}i_{A}^{(1)}(\Theta_{x^{e},\delta^{e}})\\
= & \sum_{e\in E}i_{\mathcal{X}}(x^{e})i_{\mathcal{X}}(\delta^{e})^{*}=\sum_{e\in E}\left(\sum_{f\in E}S_{f}i_{A}(x_{f}^{e})\right)\left(\sum_{g\in E}S_{g}i_{A}(\delta_{g}^{e})\right)^{*}\\
= & \sum_{e\in E}(S_{e}i_{A}(a_{s(e)}\circ\phi_{e}))(i_{A}(1_{T_{r(e)}})S_{e}^{*})=i_{A}(a)\sum_{e\in E}S_{e}S_{e}^{*}=i_{A}(a).\end{align*}
 Therefore $(i_{A},i_{\mathcal{X}})$ is a Cuntz-Pimsner covariant
representation.

For $\delta^{e}$ defined as above, we notice that\[
i_{\mathcal{X}}(\delta^{e})=\sum_{f\in E}S_{f}i_{A}(\delta_{f}^{e})=S_{e}i_{A}(1_{T_{r(e)}})=S_{e}.\]
 Then $i_{\mathcal{X}}(\mathcal{X})\cup i_{A}(A)$ generates $C^{\ast}(G)$.

Since $(i_{A},i_{\mathcal{X}})$ is a Cuntz-Pimsner covariant representation,
there exists a homomorphism $i_{\mathcal{X}}\times i_{A}$ of $\mathcal{O}_{\mathcal{X}}$
onto $C^{\ast}(G)$ such that $(i_{\mathcal{X}}\times i_{A})\circ k_{\mathcal{X}}=i_{\mathcal{X}}$
and $(i_{\mathcal{X}}\times i_{A})\circ k_{A}=i_{A}$. We prove that
$i_{\mathcal{X}}\times i_{A}$ is also injective. Let $\gamma:\mathbb{T}\to\Aut(\mathcal{O}_{\mathcal{X}})$
defined by $\gamma_{z}(k_{\mathcal{X}}(\xi))=zk_{\mathcal{X}}(\xi)$
and $\gamma_{z}(k_{A}(\xi))=k_{A}(\xi)$ be the gauge action on $\mathcal{O}_{\mathcal{X}}$.
Let $\beta:\mathbb{T}\to\Aut(C^{\ast}(G))$ defined by $\beta_{z}(S_{e})=zS_{e}$
for all $e\in E$, be the gauge action on $C^{\ast}(G)$. Therefore,
by the definition of $i_{A}$ and $i_{\mathcal{X}}$, $\beta_{z}(i_{\mathcal{X}}(\xi))=zi_{\mathcal{X}}(\xi)$
and $\beta_{z}(i_{A}(a))=i_{A}(a)$ for all $\xi\in\mathcal{X}$ and
$a\in A$. Hence $\beta_{z}\circ(i_{\mathcal{X}}\times i_{A})=(i_{\mathcal{X}}\times i_{A})\circ\gamma_{z}$
for all $z\in\mathbb{T}$. Then the Gauge-Invariant Uniqueness Theorem%
\footnote{We are grateful to the referee for suggesting the use of the Gauge-Invariant
Uniqueness Theorem here. It simplifies our original proof.%
} (\cite[Theorem 4.1]{Fowler-Muhly-Raeburn}) implies that $i_{\mathcal{X}}\times i_{A}$
is injective. Thus $C^{*}(G)$ is isomorphic to the Cuntz-Pimsner
algebra associated to the $C^{*}$-corres\-pondence $\mathcal{X}$
.\hfill{}\qed 

\begin{cor}
The Cuntz-Pimsner algebra $\mathcal{O}_{\mathcal{X}}=\mathcal{O}(J(\mathcal{X}),\mathcal{X})$
of the $C^{*}$-correspon\-dence associated with an iterated function
system $(\phi_{1},\phi_{2},...,\phi_{n})$ is isomorphic to the Cuntz
algebra $\mathcal{O}_{n}$.
\end{cor}
If $K$ (the invariant set of the Mauldin-Williams graph) is a proper
subset of $T$ then $U:=T\setminus K$ is a nonempty open set of $T$.
Let $I_{U}:=C_{0}(U)$ be the corresponding ideal in $A$. Then\[
\mathcal{X}_{I_{U}}:=\{\xi\in\mathcal{X}\;:\;\langle\xi,\eta\rangle_{A}\in I_{U}\;\mbox{for all}\;\eta\in\mathcal{X}\}\]
 is a right Hilbert $I_{U}$-module and we know that $\mathcal{X}_{I_{U}}=\mathcal{X}I_{U}:=\{\xi\cdot i\;:\;\xi\in\mathcal{X},i\in I_{U}\}$
(see \cite[Section 2]{Fowler-Muhly-Raeburn}). It follows that $\mathcal{X}_{I_{U}}=\{\xi\in\mathcal{X}\;:\;\xi_{e}\in C_{0}(U)\}$
($\xi_{e}\in C_{0}(U)$ means that $\xi(e,x)=0$ if $x\in K$). We
claim that $I_{U}$ is an $\mathcal{X}$-invariant ideal in $A$,
i.e. $\Phi(I_{U})\mathcal{X}\subset\mathcal{X}I_{U}$. For $i\in I_{U}$
and $\xi\in\mathcal{X}$ we have $(\Phi(i)\xi)_{e}=i\circ\phi_{e}\xi_{e}$,
and since $i\in I_{U}$ and $\phi_{e}(K_{r(e)})\subset K_{s(e)}$,
$i\circ\phi_{e}\in I_{U}$. Hence $(\Phi(i)\xi)_{e}\in I_{U}$. Therefore
$I_{U}$ is an $\mathcal{X}$-invariant ideal in $A$ and $\mathcal{X}/\mathcal{X}I_{U}$
is a $C^{\ast}$-correspondence over $A/I_{U}\simeq C(K)$ (see \cite[Lemma 2.3]{Fowler-Muhly-Raeburn}).
Moreover $\mathcal{X}/\mathcal{X}I_{U}\simeq\mathcal{X}(K)$, where
$\mathcal{X}(K)=C(E\times_{G}K)$ is the $C^{\ast}$-correspondence
defined as in Definition \ref{def:Mauldinwilliams module} for the
$C^{\ast}$-algebra $C(K)$. Then the ideal $\mathcal{I}(I_{U})$
of $\mathcal{O_{X}}$ generated by $i_{A}(I_{U})$ is Morita equivalent
to $\mathcal{O}_{\mathcal{X}I_{U}}$, and since $\Phi(A)\subset\mathcal{K}(\mathcal{X})$,
$\mathcal{O}_{\mathcal{X}}/\mathcal{I}(I_{U})\cong\mathcal{O}_{\mathcal{X}/\mathcal{X}I_{U}}$
(see \cite[Corollary 3.3]{Fowler-Muhly-Raeburn}).

\begin{prop}
The ideal $\mathcal{I}(I_{U})$ generated by $i_{A}(I_{U})$ is equal
to 0.
\end{prop}
\begin{proof}
Let $a\in I_{U}$ be a Lipschitz function. As in the proof of Theorem
\ref{thm:CuntzKrieger} we have that $i_{A}(a)=\lim_{k\rightarrow\infty}\sum_{\alpha\in E^{k}}\mu_{r(\alpha)}^{0}(a_{s(\alpha)}\circ\phi_{\alpha})S_{\alpha}S_{\alpha}^{*}$.
Let $\varepsilon>0$. Since $\bigcap_{k\in\mathbb{N}}\phi_{\alpha_{1}\dots\alpha_{k}}(T_{r(\alpha_{k})})=\{ x_{\alpha}\}$
with $x_{\alpha}\in K_{s(\alpha)}$ and $\diam(\phi_{\alpha_{1}\dots\alpha_{k}}(T_{r(\alpha_{k})})<c^{k}D$
for all $\alpha\in E^{\infty}$, there exists $N\in\mathbb{N}$ such
that $|\mu_{r(\alpha_{k})}^{0}(a\circ\phi_{\alpha_{1}\dots\alpha_{k}})-a(x_{\alpha})|<\varepsilon$
for all $\alpha\in E^{\infty}$ and $k\ge N$. Since $a(x)=0$ for
all $x\in K$, $\Vert\sum_{\alpha\in E^{k}}\mu_{r(\alpha)}^{0}(a_{s(\alpha)}\circ\phi_{\alpha})S_{\alpha}S_{\alpha}^{*}\Vert<\varepsilon$
for all $k\ge N$. Hence $i_{A}(a)=0$. 
\end{proof}
\begin{cor}
The Cuntz-Pimsner algebra associated to the $C^{*}$-correspondence
$C(E\times_{G}K)$ over $C(K)$ with the actions defined as in Definition
\ref{def:Mauldinwilliams module} is isomorphic to $C^{*}(G)$. 
\end{cor}
One can interpret the previous results in the particular case of the
iterated function system and obtain the result from \cite[Remark 4.6]{Pinzari-Watatami-Y}.

\section{On Noncommutative Mauldin-Williams graphs}

We give a generalization of the work of Pinzari, Watatani and Yonetani
from \cite[Section 4.3]{Pinzari-Watatami-Y} on noncommutative iterated
function systems in the context of {}``noncommutative'' Mauldin-Williams
graphs and the Rieffel metric. We show that, in fact, these situations
are no more general than those just discussed.

We begin by reviewing the Rieffel metric.

\begin{defn}
\label{def:(The-Rieffel-metric)} Let $A$ be a unital $C^{\ast}$-algebra,
let $\mathcal{L}(A)\subset A$ be a dense subspace of $A$ (the Lipschitz
elements), and let $L$ be a semi-norm (the Lipschitz seminorm) on
$\mathcal{L}(A)$ such that $\mathcal{K}:=\{ a\in\mathcal{L}(A)\;:\; L(a)=0\}$
equals the scalar multiples of the identity. The \emph{Rieffel metric}
$\rho$ on the state space $\mathcal{S}$ of $A$ is defined by the
equation\[
\rho(\mu,\nu)=\sup\{|\mu(a)-\nu(a)|\;:\; a\in\mathcal{L}(A),L(a)\leq1\}\]
 for all $\mu,\nu\in\mathcal{S}$. We will suppose that the metric
$\rho$ is bounded on $S$ and that the corresponding topology coincides
with the weak-$\ast$ topology on $\mathcal{S}$.
\end{defn}
For a compact metric space $(X,\rho)$ let\[
\mathcal{C}(X):=\{ E\;:\: E\textrm{ is a non-empty compact subset of }X\}.\]
The \emph{Hausdorff metric} on $\mathcal{C}(X)$ is defined by the
formula\[
\delta_{\rho}(E,F)=\inf\{ r>0\;:\; U_{r}(E)\supseteq F\mbox{ and }U_{r}(F)\supseteq E\}\]
for all $E,F\in\mathcal{C}(X)$, where $U_{r}(E)=\{ x\in X\;:\;\rho(x,y)<r\mbox{ for some }y\in E\}$
(see \cite[Theorem 2.4.1]{Edgar} or \cite[Proposition 1.1.5]{kigami}).
Then $(\mathcal{C}(X),\delta_{\rho})$ is a compact metric space.

For a $C^{\ast}$-algebra $A$, Rieffel defines the \emph{quantum
closed subsets} of $A$ in \cite[Page 14]{Rieffel3} to be the closed
convex subsets of the state space $\mathcal{S}(A)$ of $A$. If $L$
is a Lipschitz seminorm on $A$ and $\rho_{L}$ is the corresponding
Rieffel metric, the space $\mathcal{Q}(A)$ of quantum closed subsets
of $A$ is a compact metric space for the associated Hausdorff metric
(see \cite[Page 14]{Rieffel3}).

Following the definition of the classical Mauldin-Williams graphs,
we define a noncommutative variant.

\begin{defn}
A \emph{noncommutative Mauldin-Williams graph} is a system $(G,\{ A_{v},\mathcal{L}_{v},\allowbreak L_{v},\rho_{v}\}_{v\in V},\{\phi_{e}\}_{e\in E})$
where $G=(V,E,s,r)$ is a graph and where $\{ A_{v},\mathcal{L}_{v},L_{v}\}_{v\in V}$
and $\{\phi_{e}\}_{e\in E}$ are families such that 
\end{defn}
\begin{enumerate}
\item For each $v\in V$, $A_{v}$ is a unital $C^{\ast}$-algebra with
a prescribed Lipschitz seminorm $L_{v}$ on a prescribed subspace
$\mathcal{L}_{v}$ of Lipschitz elements in $A_{v}$ and $\rho_{v}$
is the corresponding Rieffel metric.
\item For $e\in E$, $\phi_{e}$ is a unital $\ast$-homomorphism from $A_{s(e)}$
to $A_{r(e)}$ such that\[
\rho_{s(e)}(\phi_{e}^{*}(\mu),\phi_{e}^{*}(\nu))\leq c\rho_{r(e)}(\mu,\nu)\]
 for some constant $c$ satisfying $0<c<1$ and all $\mu,\nu\in\mathcal{S}_{r(e)}$
(where $\mathcal{S}_{v}$ is the state space of the $C^{*}$-algebra
$A_{v}$). 
\end{enumerate}
We shall assume, too, that there are no sinks in the graph $G$. We
also let $\mathcal{S}:=\prod_{v\in V}\mathcal{S}_{v}$.

When we have one vertex and $n$ edges we call the system a \emph{noncommutative
iterated function system}.

Let $\mathcal{C}(\mathcal{S}_{v})$ be the space of compact subsets
of $\mathcal{S}_{v}$ endowed with the Hausdorff metric $\delta_{\rho_{L_{v}}}$,
for each $v\in V$. Let $\mathcal{C}(\mathcal{S})=\prod_{v\in V}\mathcal{C}(\mathcal{S}_{v})$.
Then $\mathcal{C}(\mathcal{S})$ is a compact metric space. Moreover,
the map $F:\mathcal{C}(\mathcal{S})\rightarrow\mathcal{C}(\mathcal{S})$
defined by the formula 
\newcommand{\twolines}[2]{{\genfrac{}{}{0pt}{}{#1}{#2}}}
\[
F\left(\left(K_{v}\right)_{v\in V}\right)=\left(\bigcup_{\twolines{e\in E}{s(e)=v}}\phi_{e}^{\ast}(K_{r(e)})\right)_{v\in V}\]
is a contraction, since each $\phi_{e}^{\ast}$ is a contraction with
respect to the Rieffel metric. Thus there exist an \emph{unique} element
$(K_{v})_{v\in V}\in\mathcal{C}(\mathcal{S})$ such that\begin{equation}
K_{v}=\bigcup_{\twolines{e\in E}{s(e)=v}}\phi_{e}^{\ast}(K_{r(e)})\label{eq:noncomminvset}\end{equation}
for all $v\in V$. Let $T_{v}$ be the closed convex hull of $K_{v}$,
for $v\in V$. That is $T_{v}\in\mathcal{Q}(A_{v})$ for all $v\in V$.
Since by \cite[Proposition 3.6]{Rieffel3} there is a bijection between
isomorphism classes of quotients of $A_{v}$ and closed convex subsets
of $\mathcal{S}_{v}$, we will assume that\begin{equation}
\mathcal{S}_{v}=T_{v}\quad\mbox{for all }v\in V,\label{eq:noncomminv2}\end{equation}
by taking a quotient of the original $C^{\ast}$-algebra $A_{v}$,
if necessary. In particular, if $(M_{v})_{v\in V}\in\mathcal{S}$
is any family which satisfies Equation (\ref{eq:noncomminvset}),
then $M_{v}=K_{v}$ and the closed convex hull of $M_{v}$ equals
$\mathcal{S}_{v}$ for all $v\in V$. 

\begin{lem}
\label{lem:invariantideal}In the above situation, if $\mathcal{I}$
is an ideal in $\sum_{v\in V}^{\oplus}A_{v}$ of the form $\mathcal{I}=(\mathcal{I}_{v})_{v\in V}$,
with $\mathcal{I}_{v}$ a proper ideal of $A_{v}$, then\begin{equation}
\mathcal{I}_{v}=\bigcap_{\genfrac{}{}{0pt}{}{e\in E}{s(e)=v}}\phi_{e}^{-1}(\mathcal{I}_{r(e)})\mbox{ if and only if }\mathcal{I}_{v}=(0_{v})\;\mbox{for all }v\in V.\label{eq:noncomminvar3}\end{equation}

\end{lem}
\begin{proof}
Let $\mathcal{I}=(\mathcal{I}_{v})_{v\in V}$ be such that $\mathcal{I}_{v}=\bigcap_{\genfrac{}{}{0pt}{}{e\in E}{s(e)=v}}\phi_{e}^{-1}(\mathcal{I}_{r(e)})$
for all $v\in V$. Let $M_{v}:=\{\mu\in\mathcal{S}_{v}\;:\;\mu(a)=0\;\mbox{for all}\; a\in\mathcal{I}_{v}\}$.
We show that $(M_{v})_{v\in V}$ is a family which satisfies Equation
(\ref{eq:noncomminvset}). Let $v\in V$ and let $\mu\in\bigcup_{\genfrac{}{}{0pt}{}{e\in E}{s(e)=v}}\phi_{e}^{\ast}(M_{r(e)})$.
Then there exists some $e\in E$ and $\nu\in M_{r(e)}$ such that
$\mu=\phi_{e}^{\ast}(\nu)$. Let $a\in\mathcal{I}_{v}$. Then $\phi_{e}(a)$
belongs to $\mathcal{I}_{r(e)}$. Hence $\mu(a)=\nu(\phi_{e}(a))=0$.
Therefore $\mu\in M_{v}$.

Now suppose that there is some $\mu\in M_{v}$ such that $\mu\notin\bigcup_{\genfrac{}{}{0pt}{}{e\in E}{s(e)=v}}\phi_{e}^{\ast}(M_{r(e)})$.
Hence there is some $a\in A_{v}$ such that $\mu(A)\ne0$ and $\phi_{e}^{\ast}(\nu)(a)=0$
for all $\nu\in M_{r(e)}$ and for all $e\in E$ such that $s(e)=v$.
Then $\phi_{e}(a_{v})\in\mathcal{I}_{r(e)}$ for all $e\in E$, therefore
$a_{v}\in\mathcal{I}_{v}$. Thus $\mu(a)=0$, which is a contradiction.
Then the family $(M_{v})_{v\in V}$ satisfies Equation (\ref{eq:noncomminvset}).
Therefore $\mu(a)=0$ for all $\mu\in\mathcal{S}_{v}$ and $a\in\mathcal{I}_{v}$,
hence $\mathcal{I}_{v}=0$ for all $v\in V$.

Suppose that there exists some $a\in A_{v}$ which is not zero, but
$a\in\bigcap_{\genfrac{}{}{0pt}{}{e\in E}{s(e)=v}}\Ker\phi_{e}$.
Then there is some $\mu\in\mathcal{S}_{v}$ such that $\mu(a)=0$.
Since $K_{v}=\bigcup_{\genfrac{}{}{0pt}{}{e\in E}{s(e)=v}}\phi_{e}^{\ast}(K_{r(e)})$,
there is some $e\in E$ with $s(e)=v$ and some $\nu\in S_{r(e)}$
such that $\mu=\phi_{e}^{\ast}(\nu)$. Since $\phi_{e}(a)=0$, we
obtain that $\mu(a)=\nu(\phi_{e}(a))=0$, which is a contradiction.
Hence $(0_{v})=\bigcap_{\genfrac{}{}{0pt}{}{e\in E}{s(e)=v}}\Ker\phi_{e}$.
\end{proof}
In a fashion similar with the commutative case, we can define a $C^{*}$-correspondence
associated to a noncommutative Mauldin-Williams graph.

\begin{defn}
Let $(G,\{ A_{v},\mathcal{L}_{v},\allowbreak L_{v},\rho_{v}\}_{v\in V},\{\phi_{e}\}_{e\in E})$
be a noncommutative Mauldin-Williams graph. Suppose that the underlying
graph has now sinks and no sources. Let $A:=\sum_{v\in V}^{\oplus}A_{v}$.
For an element $(\xi_{e})_{e\in E}$ of the column space $C^{n}(A)$,
where $n=|E|$, we have $\xi_{e}=\sum_{v\in V}^{\oplus}\xi_{e,v}$,
where $\xi_{e,v}\in A_{v}$ for all $v\in V$ and $e\in E$. We define\[
\mathcal{X}=\{\xi\in C^{n}(A)\;:\:\xi_{e,v}\equiv0\;\mbox{unless}\; v=r(e)\}\]
 and view it as a $C^{*}$-correspondence over $A$ via the formulae\[
(\xi\cdot a)_{e}=\xi_{e}\cdot a_{r(e)},\]
\[
(a\cdot\xi)_{e}=\phi_{e}(a_{s(e)})\cdot\xi_{e},\]
 where $a=\sum_{v\in V}^{\oplus}a_{v}$, and the $A$-valued inner
product defined by the formula \[
\langle\xi,\eta\rangle_{A}=\sum_{e\in E}\xi_{e}^{*}\eta_{e}\]
 for all $\xi,\eta\in\mathcal{X}$ and $a\in A$. Since the graph
$G$ has no sources, the $A$-valued inner product is well defined.

The left action is given by the $\ast$-homomorphism $\Phi:A\rightarrow\mathcal{L}(\mathcal{X})$
defined by the formula $\Phi(a)\xi=a\cdot\xi$. Since $\Phi(a)=0$
if and only if $a_{v}\in\bigcap_{s(e)=v}\Ker\phi_{e}$ for all $v\in V$,
Lemma \ref{lem:invariantideal} implies that the left action is faithful.
\end{defn}
In the noncommutative iterated function system case, the $C^{\ast}$-correspondence
is $\mathcal{X}=C^{n}(A)$ with left action given by the $\ast$-homomorphism
$a\rightarrow\operatorname{diag}(\phi_{i}(a))$.

Recall that $E^{k}$ denotes the set of paths of length $k$, $E^{\infty}$
denotes the set of infinite paths in the graph $G$, $E^{k}(v)$ denotes
the set of paths of length $k$ starting at the vertex $v$, and $E^{\infty}(v)$
denotes the set of infinite paths starting at the vertex $v$. For
$k\in\mathbb{N}$ and $\alpha\in E^{k}$, we write $\phi_{\alpha_{1}\dots\alpha_{k}}^{*}$
for the map $\phi_{\alpha_{1}}^{\ast}\circ\dots\circ\phi_{\alpha_{k}}^{\ast}:\mathcal{S}_{r(\alpha_{k})}\to\mathcal{S}_{s(\alpha_{1})}$
and $\phi_{\alpha_{k}\dots\alpha_{1}}$ for the map $\phi_{\alpha_{k}}\circ\dots\circ\phi_{\alpha_{1}}:A_{s(\alpha_{1})}\to A_{r(\alpha_{k})}$.
We will use the following results (which are similar to the commutative
case): if $v\in V$, $a\in\mathcal{L}_{v}(A_{v})$ and $\mu,\nu\in\mathcal{S}_{v}$
then $|\mu(a)-\nu(a)|\leq\rho_{v}(\mu,\nu)\cdot L_{v}(a)$; if $\alpha\in E^{k}$
and $\mu,\nu\in\mathcal{S}_{r(\alpha)}$, then\begin{equation}
\rho_{s(\alpha)}(\phi_{\alpha}^{\ast}(\mu),\phi_{\alpha}^{\ast}(\nu))\leq c^{k}\rho_{r(\alpha)}(\mu,\nu)\leq c^{k}D,\label{eq:statesineq}\end{equation}
 where $s(\alpha)=s(\alpha_{1})$, $r(\alpha)=r(\alpha_{k})$ and
$D=\max_{v\in V}\operatorname{diam}_{L_{v}}(\mathcal{S}_{v})$.

Since $(G,\{\mathcal{S}_{v},\rho_{v}\}_{v\in V},\{\phi_{e}^{\ast}\}_{e\in E})$
is a (classical) Mauldin-Williams graph, for each $\alpha\in E^{\infty}$
there is a unique state $\mu_{\alpha}\in\mathcal{S}_{s(\alpha)}$
such that $\{\mu_{\alpha}\}=\bigcap_{k\in\mathbb{N}}\phi_{\alpha_{1}\dots\alpha_{k}}^{\ast}(\mathcal{S}_{r(\alpha_{k})})$.
In particular $\lim_{k\rightarrow\infty}\phi_{\alpha_{1}\cdots\alpha_{k}}^{*}(\mu_{r(\alpha_{k})})=\mu_{\alpha}$
for all $\mu=(\mu_{v})_{v\in V}\in\mathcal{S}$.

\begin{thm}
\label{thm:homintocomm}Let $(G,\{ A_{v},\mathcal{L}_{v},L_{v},\rho_{v}\}_{v\in V},\{\phi_{e}\}_{e\in E})$
be a (noncommutative) Mauldin-Williams graph. Suppose that the graph
$G$ has no sinks. Then there is an injective $\ast$-homomorphism
from $A$ into $C(E^{\infty})$.
\end{thm}
\begin{proof}
Fix $v_{0}\in V$. Define $\pi_{v_{0}}:A_{v_{0}}\rightarrow C(E^{\infty}(v_{0}))$
by the formula\[
\pi_{v_{0}}(a)(\alpha)=\mu_{\alpha}(a)\]
for all $a\in A_{v_{0}}$. Thus, if $a\in\mathcal{L}_{v_{0}}$,\[
\pi_{v_{0}}(a)(\alpha)=\lim_{k\to\infty}\mu_{r(\alpha_{k})}(\phi_{\alpha_{k}\cdots\alpha_{1}}(a))\]
for all $\mu=(\mu_{v})_{v\in V}\in\mathcal{S}$. By the comments preceding
the theorem the map $\pi_{v_{0}}$ is well defined. We prove that
it is a homomorphism.

Let $\mu_{0}=(\mu_{v}^{0})_{v\in V}\in\mathcal{S}$ be \emph{fixed}.
Let $a\in\mathcal{L}_{v_{0}}$, $\alpha\in E^{\infty}(v_{0})$ and
let $\varepsilon>0$. Let $k\in\mathbb{N}$ such that $c^{k}DL_{v_{0}}(a)<\varepsilon$.
For \emph{any} $\mu=(\mu_{v})_{v\in V}\in\mathcal{S}$ we have\begin{alignat*}{1}
\left|\mu_{r(\alpha_{k})}\left(\mu_{r(\alpha_{k})}^{0}(\phi_{\alpha_{k}\cdots\alpha_{1}}(a))1_{A_{r(\alpha_{k})}}-\phi_{\alpha_{k}\cdots\alpha_{1}}(a)\right)\right|\end{alignat*}
\[
=\left|\mu_{r(\alpha_{k})}^{0}(\phi_{\alpha_{k}\cdots\alpha_{1}}(a))-\mu_{r(\alpha_{k})}(\phi_{\alpha_{k}\cdots\alpha_{1}}(a))\right|<c^{k}DL_{v_{0}}(a)<\varepsilon.\]
Hence \begin{equation}
\Vert\mu_{r(\alpha_{k})}^{0}(\phi_{\alpha_{k}\cdots\alpha_{1}}(a))1_{A_{r(\alpha_{k})}}-\phi_{\alpha_{k}\cdots\alpha_{1}}(a)\Vert<4\varepsilon.\label{eq:ineq1}\end{equation}
Let $a,b\in A_{v_{0}}$. We have\[
|\mu_{r(\alpha_{k})}^{0}(\phi_{\alpha_{k}\cdots\alpha_{1}}(ab))-\mu_{r(\alpha_{k})}^{0}(\phi_{\alpha_{k}\cdots\alpha_{1}}(a))\mu_{r(\alpha_{k})}^{0}(\phi_{\alpha_{k}\cdots\alpha_{1}}(b))|\]
\begin{eqnarray*}
 & \le & \Vert\mu_{r(\alpha_{k})}^{0}(\phi_{\alpha_{k}\cdots\alpha_{1}}(ab))-\phi_{\alpha_{k}\cdots\alpha_{1}}(a)\phi_{\alpha_{k}\cdots\alpha_{1}}(b)\Vert\\
 & + & \Vert\phi_{\alpha_{k}\cdots\alpha_{1}}(a)\phi_{\alpha_{k}\cdots\alpha_{1}}(b)-\phi_{\alpha_{k}\cdots\alpha_{1}}(a)\mu_{r(\alpha_{k})}^{0}(\phi_{\alpha_{k}\cdots\alpha_{1}}(b))\Vert\\
 & + & \Vert\phi_{\alpha_{k}\cdots\alpha_{1}}(a)\mu_{r(\alpha_{k})}^{0}(\phi_{\alpha_{k}\cdots\alpha_{1}}(b))-\mu_{r(\alpha_{k})}^{0}(\phi_{\alpha_{k}\cdots\alpha_{1}}(a))\mu_{r(\alpha_{k})}^{0}(\phi_{\alpha_{k}\cdots\alpha_{1}}(b))\Vert\\
 & < & 4\varepsilon+4\Vert a\Vert\varepsilon+4\varepsilon=(8+4\Vert a\Vert)\varepsilon,\end{eqnarray*}
by Inequality \ref{eq:ineq1}. Since $\lim_{k\to\infty}\mu_{r(\alpha_{k})}^{0}(\phi_{\alpha_{k}\cdots\alpha_{1}}(ab))=\pi_{v_{0}}(ab)(\alpha)$
and\newline $\lim_{k\to\infty}\mu_{r(\alpha_{k})}^{0}(\phi_{\alpha_{k}\cdots\alpha_{1}}(a))\mu_{r(\alpha_{k})}^{0}(\phi_{\alpha_{k}\cdots\alpha_{1}}(b))=\pi_{v_{0}}(a)(\alpha)\pi_{v_{0}}(b)(\alpha)$,
we see that $\pi_{v_{0}}$ is a homomorphism.

Hence, for each $v\in V$, we have defined a $\ast$-homomorphism
$\pi_{v}:A_{v}\rightarrow C(E^{\infty}(v))$. We prove that $\pi_{v}$
is injective for all $v\in V$. Let $v\in V$. Let $a\in A_{v}$.
Then\begin{eqnarray*}
a\in\Ker\pi_{v} & \;\Leftrightarrow & \pi_{v}(a)(\alpha)=0\;\mbox{for all }\alpha=(\alpha_{n})_{n\in\mathbb{N}}\in E^{\infty}(v)\\
 & \Leftrightarrow & \pi_{v}(a)(e\beta)=0\;\mbox{for all }\beta\in E^{\infty}(r(e))\mbox{ and }e\in E(v)\\
 & \Leftrightarrow & \pi_{r(e)}(\phi_{e}(a))(\beta)=0\;\mbox{for all }\beta\in E^{\infty}(r(e))\mbox{ and }e\in E(v)\\
 & \Leftrightarrow & a\in\bigcap_{s(e)=v}\Ker\pi_{r(e)}\circ\phi_{e}.\end{eqnarray*}
 Lemma \ref{lem:invariantideal} implies that $\Ker\pi_{v}=0$ for
all $v\in V$, hence the $\ast$-homomorphism $\pi:A\rightarrow C(E^{\infty})$
defined by the formula\[
\pi\left((a_{v})_{v\in V}\right)=\sum_{v\in V}\!^{\oplus}\pi_{v}(a_{v})\]
is an injective $\ast$-homomorphism.
\end{proof}
\begin{cor}
Under the hypothesis of Theorem \ref{thm:homintocomm}, we conclude
that $A$ must be a commutative $C^{*}$-algebra. 
\end{cor}
Even in the setting of a {}``noncommutative'' iterated function
system studied in \cite[Section 4.2]{Pinzari-Watatami-Y}, if we have
defined a Rieffel metric such that the underlying topology and the
weak-$\ast$ topology coincide and if the duals of the endomorphisms
restricted to the state space of $A$ are contractions with respect
to the Rieffel metric, then (under the hypothesis that $A$ satisfies
Equation \ref{eq:noncomminv2}) $A$ is forced to be commutative and
the endomorphisms $\phi_{i}$ must come from an ordinary iterated
function system, i.e. $A=C(K)$ for some compact metric space and
there are contractions $\{\varphi_{i}\}_{i=1,\dots,n}$ defined on
$K$ such that $\phi_{i}(a)=a\circ\varphi_{i}$. This seems not to
have been noticed by the authors of \cite{Pinzari-Watatami-Y}.

The assumption that the graph $G$ has no sinks is essential in the
proof of Theorem \ref{thm:CuntzKrieger} and Theorem \ref{thm:homintocomm},
since it forces the presence of infinite paths in the graph. Also,
the assumption that the graph $G$ has no sources was needed to define
the $C^{\ast}$-correspondence associated with a Mauldin-Williams
graph. It is not needed, though, in the proof of Theorem \ref{thm:homintocomm}.

We would like to call attention to a recent preprint of Kajiwara and
Watatani \cite{KWp03} in which they considered a somewhat different
$C^{\ast}$-correspondences associated with an iterated function system
and arrive at a $C^{\ast}$-algebra that is sometimes different from
$\mathcal{O}_{n}$. It appears that their construction can be modified
to cover the setting of Mauldin-Williams graphs, leading to $C^{\ast}$-algebras
different from the Cuntz-Krieger algebras of the underlying graphs.
We intend to pursue the ramifications of this in a future note.

\end{document}